\documentclass{article}

\usepackage[latin1]{inputenc}
\usepackage{amsmath, amsthm, amsfonts}
\usepackage[english]{babel}
\usepackage{graphics}
\usepackage{graphicx}
\usepackage{amsfonts}
\usepackage{amssymb}
\usepackage{comment}
\usepackage{subfigure}
\usepackage{ragged2e}

\theoremstyle{definition}

\theoremstyle{remark}

\newtheorem{defi}{Definition}[section]

\newtheorem*{teo2}{Keywords}

\newtheorem{lema}[defi]{Lemma}
\newtheorem{teo}[defi]{Theorem}

\theoremstyle{definition}

\theoremstyle{remark}
\newtheorem*{remark}{Remark}

\title{Generalized Ribaucour-type surfaces }
\author{Milton Javier Cardenas Mendez \\ Instituto de Matemática e Estatística \\
Universidade Federal de Goiás, 74001-970 Goiânia, GO, Brazil\\ miltonjcardenas@ufg.br \and  Armardo Mauro Vasquez Corro\\ Instituto de Matemática e Estatística \\ Universidade Federal de Goiás, 74001-970 Goiânia, GO, Brazil\\corro@ufg.br }
\date{\today}

\begin{document}
\maketitle

\begin{abstract}
In this work we generalize the surfaces studied in \cite{mi}, we define the generalization of Ribaucour-type surfaces (in short, GRT-surfaces). We obtain present a representation for GRT-surfaces with prescribed Gauss map which depends on two holomorphic functions and a real function  $\ell$. We give explicit examples of GRT-surfaces. Also, we use this representation to classify the GRT-surfaces of rotation.
\end{abstract}

\begin{teo2}
Generalitation of Ribaucour-type surfaces, generalized Weingarten surfaces, prescribed normal Gauss map.
\end{teo2}

\textbf{Introduction}: Let $\Sigma\subset \mathbb{R}^3$ be an oriented surface with normal Gauss map $N$, functions $\Psi,\Lambda :\Sigma\rightarrow \mathbb{R}^3$ given by $ \Psi(p) = \langle p , N(p) \rangle$ and $ \Lambda(p) = \langle p , p \rangle $, $ p \in \Sigma$, where $\langle,\rangle$ denotes the Euclidean scalar product in $\mathbb{R}^3$, are called support function 
and quadratic distance function, respectively. Geometrically, $\Psi(p)$ measures the signed distance from the origin to $T_pM$ and $\Lambda(p)$ measures the squared distance from the origin to $p$. Let $p\in \Sigma $, a sphere with center $p+\frac{H(p)}{K(p)} N(p)$ and radius $\frac{H(p)}{K(p) }$ is called the middle sphere.

In 1888, Appell \cite{f} studied a class of surfaces oriented in $ \mathbb{R}^3 $ associated with area-preserving sphere transformations. Later, Ferreira and Roitman \cite{d} showed that these surfaces satisfy the Weingarten relation, $H + \Psi K = 0 $.

In 1907, Tzitzeica \cite{o} studied hyperbolic surfaces oriented so that there is a nonzero constant $c \in \mathbb{R} $ for which $ K + c^2 \Psi^4 = 0 $.

In \cite{a}, the authors motivated by the works \cite{f} and \cite{o} defined  generalized Weingarten surfaces as surfaces that satisfy a relation of the form $A + BH + CK = 0$, where $A, B, C : \Sigma \rightarrow \mathbb{R}$ are differentiable functions that do not depend on the parameterization of $\Sigma$.
  In particular, they studied the class of surfaces that satisfy the relation $2\Psi H + \Lambda K = 0$. Called Special Generalized Weingartem Surfaces depending on the support function and the distance function (in short, EDSGW-surfaces), these surfaces have the geometric property that all medium spheres pass through the origin. The authors obtained a Weierstrass-like representation of EDSGW-surfaces depending on two homomorphic functions. In \cite{b}, the authors classified isothermic EDSGW-surfaces in relation to the third fundamental shape parameterized by plane curvature lines. Also in \cite{z}, it is shown that EDSGW-surfaces are in correspondence with the surface class in $\mathbb{S}^2 \times \mathbb{R}$ where the Gaussian curvature $ K$ and the extrinsic curvature $K_E$ satisfy $K=2K_E$.

Martínez and Roitman, in \cite{c} showed what appears to be the first example found for the second case of the problem posed by Élie Cartan in his classic book about external differential systems and their applications to Differential Geometry. Such examples are given by a class of Weingarten surfaces that satisfy the relation $2 \Psi H + (1+ \Lambda) K = 0 $, Ribaucour surface calls, these surfaces have the geometric property that all the medial spheres intercept a fixed sphere along a large circle.

In \cite{pp}, the authors define a surface class called  Ribaucour surface of harmonic type (in short, HR-surface) if it satisfies
$2\Psi H + (c e^{2\mu} +\Lambda)K = 0$, where $c$ is a nonzero real constant, $\mu$ a harmonic function with respect to the third
fundamental form. These surfaces generalize  Ribaucour surfaces studied in \cite{c}.

Motivated by \cite{mi}, we define  $\Sigma$ be a surface with Gauss map $N$ is called a \textbf{surface of type Ribaucour generalized} or abbreviated GRT-Surface. There is a harmonic map $\mu$ with respect to the third fundamental form and a function $C:\mathbb{R}\rightarrow \mathbb{R}$ such that for all $p\in \Sigma$ the sphere with center $p+\left(\frac{H}{C(\mu)K}+\frac{\Psi(2-C(\mu))}{2C(\mu)}\right)N(p)$ and radius $\left(\frac{H}{C(\mu)K}+\frac{\Psi(2-C(\mu))}{2C(\mu)}\right)$, tangents a fixed sphere. In this case $\Sigma$ satisfies the following relation generalized Weingarten
\begin{equation*}
\frac{H}{K}=C(\mu) \left(\frac{-\Lambda}{2\Psi}+ \frac{\Psi}{2} \right)-\Psi.
\end{equation*}
 We obtain present a representation for GRT-surfaces with prescribed Gauss map which depends on two holomorphic functions and a real function  $\ell$. We give explicit examples of GRT-surfaces. Also, we use this representation to classify the GRT-surfaces of rotation.

\section{Preliminaries}
In this section we fix the notation used in this work, $ \Sigma $ a surface on $ \mathbb{R}^3 $, $ N $ its normal Gaussian map, and $ U $ an open subset of $ \mathbb{R}^2 $.

Let $ X: U \subset \mathbb{R}^2 \rightarrow \Sigma$, a parameterization of a surface $ \Sigma $ and $ N: U \subset \mathbb{R}^2 \rightarrow \mathbb{R}^{3} $, normal Gaussian map. Considering $ \{X_{,1}, X_{,2}, N \} $ as a base of $ \mathbb{R}^{3} $,  where $ X_{,i} (q) = \frac{\partial X} {\partial u_i} (q) $, $ 1 \leq i, j \leq 2 $ further we can write  vector $ X_{, ij} $, $1 \leq i, j \leq 2 $, as
\begin{equation*}
  X_{,ij}= \sum_ {k = 1}^ {2} \widetilde{\Gamma}_ {ij}^k X_{,k} + b_{ij} N
\end{equation*}
The $ \widetilde{\Gamma}_ {ij}^k $ coefficients are called Christoffel symbols.

\begin{defi}
\label{def1}
Let $ X $ be a local parameterization of $ \Sigma $ with map of Gauss $ N $, matrix $ W = (W_{ij})$, such that
\begin{equation*}
N_{,i} = \sum_{j = 1}^{2} W_{ij} X_j, \hspace{0.3cm} 1 \leq i \leq 2
\end{equation*}
is called the Weingarten matrix of $\Sigma$.
\end{defi}

\begin{lema}
\label{le1}
Let $N$ be the normal Gaussian map given by (\ref{N1})  such that the metric $ L_{ij} =\langle N_{, i}, N_{, j} \rangle $ is Euclidean conformal.  Christoffel's symbols for metric $ L_ {ij}$ are given by
\begin{equation*}
\Gamma_{ij}^k=0, \hspace{0.3cm} \Gamma_{ii}^i=\frac{L_{ii,i}}{2L_{ii}},\hspace{0.3cm} \Gamma_{ij}^i=\frac{L_{ii,j}}{2L_{ii}},\hspace{0.3cm} \Gamma_{ii}^i=\frac{-L_{ii,j}}{2L_{jj}}
\end{equation*}
For $i,j$ e $k$  different.
\end{lema}

\begin{remark}
\label{holo}
 Let the inner product be defined by $ \langle, \rangle:\mathbb{C} \times \mathbb{C} \rightarrow \mathbb{R}, \langle f,g \rangle= f_1g_1+f_2g_2$, where $f= f_1+if_2$ and $g=g_1+ig_2$ are holomorphic functions. If $f,g,h:U\subset\mathbb{C}\rightarrow\mathbb{C}$ are holomorphic functions, then
\begin{equation}\label{holo}
     \begin{split}
         &\langle f,g\rangle_{,1}=\langle f',g\rangle+\langle f,g'\rangle  \\
                    &\langle f,g\rangle_{,2}=\langle if',g\rangle+\langle f,ig'\rangle\\
             &\langle fh,g\rangle=\langle f,\overline{h} g\rangle\\
    &f =\langle 1,g\rangle+i\langle i,f\rangle\\
 & g_{,1}=g', g_{,2}=ig'\\
  \end{split}
\end{equation}
where  $\langle f,g\rangle_{,1}=\frac{\partial\langle f(z),g(z)\rangle}{\partial u_1} $
 \end{remark}

Below we present some important results studied in \cite{mi}. 
\begin{teo}
\label{te2}
Let $\Sigma \subset \mathbb{R}^3 $, an orientable surface with non-zero Gauss-Kronecker  curvature. Then there is  a differentiable  function  $h:U \rightarrow \mathbb{R}$ and  $g$ a holomorphic function  such that normal Gauss map  is given by
 \begin{equation}
\label{N1}
  N=\frac{( 2g(u),1-|g(u)|^2)}{1+|g(u)|^2}
\end{equation}
the coefficients of the $III$ fundamental form are
\begin{equation}
\label{simb}
L_{ij}=\frac{4 |g'|^ 2 \delta _{ij} }{(1+|g|^ 2)^ 2}
\end{equation}
 $\Sigma $ is locally parameterized by
\begin{equation}
\label{1}
X(u)=\sum_{j=1}^{2} \frac{h(u)_{,j}}{L_{jj}}N(u)_{,j}+h(u)N(u)
\end{equation}
In this case   $h(u)=\langle X(u),N(u)\rangle$ is the support function. Furthermore, the Weingarten matrix is given by
$ W = V^{-1} $ where
\begin{equation*}
\label{2}
V_{ij}=\frac{1}{L_{ij}}\left( h_{,ij}-\sum_{k}^{n}h_{,k}\Gamma_{ij}^{k}+hL_{ij}\delta_{ij}\right)
\end{equation*}
where $\Gamma_{ij}^{k}$ are Christoffel's symbols  of $N$ and the fundamental forms $ I $ and $ II $ of $ X $, in local coordinates, are given by

 \begin{equation*}
 I=\langle X_{,i},X_{,j}\rangle=\sum_{k=1}^{n} V_{ik}V_{jk}L_{kk},\quad II=\langle X_{,i},N_{,j}\rangle=V_{ij}L_{jj}.
 \end{equation*}
 furthermore, 
\begin{equation}
\label{lapla}
 \frac{\triangle h}{L_{11}}+2h=\frac{-2H}{K},
\end{equation}
\begin{equation}
\label{qu}
\Lambda(q)=\langle X(q), X(q)\rangle=|\nabla_Lh(q)|^2+h(q)^2
\end{equation}
\begin{equation*}
\Psi(q)=\langle X(q), N(q)\rangle=h(q)
\end{equation*}
called quadratic distance function  and support function, respectively.
\end{teo}

\section{GRT-surfaces}
In this chapter, we introduce a few classes of surfaces called surfaces.
of generalized Ribaucour type and we call it GRT-surfaces, we show a
local parameterization of this class of surfaces and characterize the case where such
surfaces are rotating.

\begin{defi}
Let $\Sigma$ be a surface with Gauss map $N$ is called a \textbf{surface of type Ribaucour generalized} or abbreviated GRT-Surface. There is a harmonic map $\mu$ with respect to the third fundamental form and a function $C:\mathbb{R}\rightarrow \mathbb{R}$ such that for all $p\in \Sigma$ the sphere with center $p+\left(\frac{H}{C(\mu)K}+\frac{\Psi(2-C(\mu))}{2C(\mu)}\right)N(p)$ and radius $\left(\frac{H}{C(\mu)K}+\frac{\Psi(2-C(\mu))}{2C(\mu)}\right)$, tangents a fixed sphere. In this case $\Sigma$ satisfies the following relation generalized Weingarten
\begin{equation}
\label{m1}
\frac{H}{K}=C(\mu) \left(\frac{-\Lambda}{2\Psi}+ \frac{\Psi}{2} \right)-\Psi.
\end{equation}
\end{defi}

\begin{remark}
We show particular cases of GRT-surfaces for different values of $C$.
\begin{enumerate}
\item If $C(t)=0$ then 
  \begin{equation*}
  \frac{H}{K}=C(t)\left(\frac{-\Lambda}{2\Psi}+ \frac{\Psi}{2} \right)-\Psi=-\Psi
\end{equation*}
  \begin{equation*}
   H+ \Psi K=0.
\end{equation*}
These are the Appell surfaces associated with area-preserving transformations on the sphere, studied by  Ferreira and Roitman  in \cite{d}.
\item If $C(t)=1$ then
  \begin{equation*}
\frac{H}{K}=C(t)\left(\frac{-\Lambda}{2\Psi}+ \frac{\Psi}{2} \right)-\Psi=\frac{-\Lambda}{2\Psi}+ \frac{\Psi}{2}-\Psi=-\frac{\Lambda}{2\Psi}- \frac{\Psi}{2}.
\end{equation*}
They are TR-surfaces studied in \cite{mi}	
\end{enumerate}
\end{remark}

\begin{lema}
Consider holomorphic functions $g:\mathbb{C}\rightarrow\mathbb{C}_{\infty}$ and $f:\Sigma \rightarrow\mathbb{C}$, with $g'\neq 0$, where $ \Sigma$ is a Riemann surface. Taking the local parameters $z=u_1+iu_2 \in \Sigma$ and $\mu=\langle 1,f\rangle$ such that $h=\ell(\mu)$, the matrix
\begin{equation*}
V_{ij}=\frac{1}{L_{jj}}\left(h_{,ij}-\sum_{k=1}^{2}h_{,k}\Gamma_{ij}^k+hL_ {ij}\delta_{ij}\right), \hspace{0.4cm} 1\leq i,j\leq 2.
\end{equation*}
Using the metric given by (\ref{simb}), we have
\begin{equation}
\label{vs}
\begin{split}
    &V_{11} =\frac{T^2}{4|g'|^2}[\ell''(\mu) \langle 1,f'\rangle^2- \ell'(\mu)\langle 1,\xi \rangle] + \ell(\mu), \\
      &V_{12}=V_{21}=\frac{T^2}{4|g'|^2}\left[\ell''(\mu) \left\langle 1,\frac{if'^2 }{2}\right\rangle + \ell'(\mu)\langle i,\xi \rangle\right],\\
      &V_{22}=\frac{T^2}{4|g'|^2}[\ell''(\mu) \langle 1,if'\rangle^2+ \ell'(\mu)\langle 1,\xi \rangle] + \ell(\mu).
\end{split}
\end{equation}
Here $\xi =f' \left(\frac{g''}{g'}-\frac{2}{T}g'\overline{g}\right)-f''$, furthermore, trace of $V$ give for 
\begin{equation*}
       trV=\frac{\ell''(\mu)|f'|^2T^2}{4|g'|^2} + 2\ell(\mu).
\end{equation*}
\end{lema}

\begin{proof}
Let $h=\ell(\mu)=\ell(\langle 1,f\rangle)$, the derivatives of $h$ are given by:
\begin{equation*}
\begin{split}
&h_{,1}=\ell'(\mu) \langle 1,f'\rangle\\
    &h_{,11} = \ell''(\mu) \langle 1,f'\rangle^2 +\ell'(\mu) \langle 1,f''\rangle \\
    &h_{,2}=\ell'(\mu) \langle 1,if'\rangle\\
      &h_{,22} =\ell''(\mu) \langle 1,if'\rangle^2 -\ell'(\mu) \langle 1,f''\rangle\\
      &h_{,12} = \ell''(\mu)\langle 1,if'\rangle\langle 1,f'\rangle+ \ell'(\mu) \langle 1,if''\rangle.
\end{split}
\end{equation*}

Using the Christoffel symbols found in Lemma \ref{le1}, we can infer that
\begin{equation*}
    \Gamma_{11}^1 =\frac{T\langle g',g''\rangle-2|g'|^2\langle g,g'\rangle }{T|g'|^2},
\end{equation*}
\begin{equation*}
\Gamma_{22}^2=\frac{T\langle g',ig''\rangle-2|g'|^2\langle g,ig'\rangle }{T|g'|^2},
\end{equation*}
\begin{equation*}
\Gamma_{11}^2=\frac{2|g'|^2\langle g,ig'\rangle -T\langle g',ig''\rangle}{T|g'|^2},
\end{equation*}
\begin{equation*}
\Gamma_{22}^1=\frac{2|g'|^2\langle g,g'\rangle-T\langle g',g''\rangle }{T|g'|^2}.
\end{equation*}

Using the previous expressions and  (\ref{holo}) we have
\begin{equation*}
\begin{split}
    \sum_{k=1}^{2}\Gamma_{11}^k h_{,k} &=\Gamma_{11}^1 h_{,1}+\Gamma_{11}^2 h_{,2} \\
    &=\ell'(\mu)\langle 1,f'\rangle \left[\frac{T\langle g',g''\rangle-2|g'|^2\langle g,g'\rangle }{T|g'|^2}\right]\\
    &-\ell'(\mu)\langle i,f'\rangle\left[\frac{2|g'|^2\langle g,ig'\rangle -T\langle g',ig''\rangle }{T|g'|^2}\right]\\
      &=\ell'(\mu) \left\langle \langle 1,f'\rangle +i\langle i,f'\rangle, \frac{\langle g',g''\rangle}{|g' |^2}-\frac{2}{T}\langle g,g'\rangle + i\left(\frac{\langle g',ig''\rangle}{|g'|^2}-\frac{2}{T}\langle g,ig'\rangle\right)\right\rangle\\
          &=\ell'(\mu)\left\langle f',\frac{\langle g',g''\rangle}{|g'|^2}+i\frac{\langle g',ig''\rangle}{|g'|^2}-\frac{2}{T}(\langle g,g'\rangle+i\langle g,ig'\rangle)\right\rangle\\
          &=\ell'(\mu)\left\langle f',\frac{\bar{g''}g'}{|g'|^2}-\frac{2}{T}\bar{g '}g\right\rangle=\ell'(\mu)\left\langle f',\overline{\left(\frac{g''}{g'}\right)}-\frac{2}{ T}\overline{g'}g\right\rangle\\
          &=\ell'(\mu)\left\langle f',\overline{\frac{g''}{g'}-\frac{2}{T}\overline{g}g'}\right\rangle=\ell'(\mu)\left\langle1,f' \left(\frac{g''}{g'}-\frac{2}{T}g'\overline{g}\right)\right\rangle\\
          &=\ell'(\mu)\langle 1,\xi+f''\rangle,
\end{split}
\end{equation*}

where $\xi =f' \left(\frac{g''}{g'}-\frac{2}{T}g'\overline{g}\right)-f''$. Similarly, knowing that $\Gamma_{21}^1=\Gamma_{22}^2, \Gamma_{12}^2=\Gamma_{11}^1$, we obtain
\begin{equation*}
\sum_{k=1}^{2}\Gamma_{22}^kh_{,k}=-\ell'(\mu)\langle 1,\xi+f'' \rangle,
\end{equation*}
It is
\begin{equation*}
\sum_{k=1}^{2}\Gamma_{12}^kh_{,k}=-\ell'(\mu)\langle i,\xi +f''\rangle.
\end{equation*}
So we can write

\begin{equation*}
\begin{split}
    V_{11} &=\frac{1}{L_{11}}\left(h_{,11}-\sum_{k=1}^{2}\Gamma_{11}^kh_{,k}+hL_ {11}\right) \\
      &=\frac{T^2}{4|g'|^2}[\ell''(\mu) \langle 1,f'\rangle^2 +\ell'(\mu) \langle 1,f ''\rangle -\ell'(\mu)\langle 1,\xi+f'' \rangle] + \ell(\mu)\\
      &=\frac{T^2}{4|g'|^2}[\ell''(\mu) \langle 1,f'\rangle^2+ \ell'(\mu)[\langle 1, f''\rangle-\langle 1,\xi+f'' \rangle] ] + \ell(\mu)\\
      &=\frac{T^2}{4|g'|^2}[\ell''(\mu) \langle 1,f'\rangle^2- \ell'(\mu)\langle 1,\ xi \rangle] + \ell(\mu)
\end{split}
\end{equation*}
\begin{equation*}
\begin{split}
V_{22} &=\frac{1}{L_{22}}\left(h_{,22}-\sum_{k=1}^{2}\Gamma_{22}^kh_{,k}+hL_ {22}\right) \\
      &=\frac{T^2}{4|g'|^2}[\ell''(\mu) \langle 1,if'\rangle^2+ \ell'(\mu)\langle 1,\ xi \rangle] + \ell(\mu)
\end{split}
\end{equation*}
\begin{equation*}
\begin{split}
V_{12} &=\frac{1}{L_{22}}\left(h_{,12}-\sum_{k=1}^{2}\Gamma_{12}^kh_{,k}\right ) \\
      &=\frac{T^2}{4|g'|^2}\left[\ell''(\mu) \left\langle 1,\frac{if'^2}{2}\right\rangle + \ell'(\mu)\langle i,\xi \rangle\right].
\end{split}
\end{equation*}
Thus,
\begin{equation*}
\begin{split}
    trV &= V_{11}+V_{22} \\
      &=\frac{T^2}{4|g'|^2}[\ell''(\mu)(\langle 1,f'\rangle^2+\langle 1,if' \rangle^2) ] + 2\ell(\mu)
\end{split}
\end{equation*}
therefore,
\begin{equation*}
trV =\frac{\ell''(\mu)|f'|^2T^2}{4|g'|^2} + 2\ell(\mu).
\end{equation*}
\end{proof}

\begin{teo}
  Let $X$ be given by (\ref{1}), $C:\mathbb{R}\rightarrow \mathbb{R}$ and $\mu:\mathbb{R}^2 \rightarrow \mathbb{R} $, $\mu$ harmonic with respect to third fundamental form then $X$ is a TRG-surface if and only if
  \begin{equation}
  \label{lapla1}
  h\triangle_Lh-C(\mu)|\nabla_Lh|^2 =0.
  \end{equation}
   \end{teo}
\begin{proof}
Using the equation (\ref{lapla}) we have
   \begin{equation}
   \label{d4}
  \frac{H}{K}=-\frac{1}{2} (\triangle_Lh+2h).
   \end{equation}
   Using (\ref{m1}) and (\ref{d4}), we get the following equivalence
   \begin{equation*}
   -\frac{1}{2} (\triangle_Lh+2h)= C(\mu) \left(\frac{-\Lambda}{2\Psi}+ \frac{\Psi}{2} \right)- \Psi=\frac{C(\mu)}{2\Psi}(-\Lambda+\Psi^2)-\Psi
   \end{equation*}

   \begin{equation*}
   \Leftrightarrow \triangle_Lh+2h= - \frac{C(\mu)}{\Psi}(-\Lambda+\Psi^2)+2\Psi,
   \end{equation*}
   remembering that $h=\Psi$, we arrive at
   \begin{equation*}
  h\triangle_Lh=C(\mu)(\Lambda-\Psi^2),
   \end{equation*}
by (\ref{qu}) we get
\begin{equation*}
  h\triangle_Lh=C(\mu)(|\nabla_Lh|^2+\Psi^2-\Psi^2),
   \end{equation*}
logo $ h\triangle_Lh-C(\mu)|\nabla_Lh|^2=0$.
  \end{proof}

The following theorem allows to obtain a class of GRT-surfaces.
  \begin{teo}
  \label{ttrg}
Let $\Sigma$ be a Riemann surface and $X:\Sigma \rightarrow \mathbb{R}^3$ be an immersion such that the Gauss-kronecker curvature is non-zero. Let $\ell:\mathbb{R}\rightarrow \mathbb{R}$ and $\mu=\langle 1,f\rangle$, with $f$ a holomorphic function and $h=\ell(\mu)$ . Then $X$ is given by (\ref{1}) and is a TRG-surface with $C=\frac{ \ell \ell''}{(\ell')^2}$, where $\ell' \neq 0$. Furthermore, let $g$ be a holomorphic function and $g'\neq 0$, then $X(\Sigma)$ is locally parameterized by
\begin{equation}
\label{trg}
X=\frac{\ell'(\mu)}{2|g'|^2}(Tg'\bar{f'}-2g\langle g',gf'\rangle , -2\langle g', gf'\rangle) + \ell(\mu)\frac{(2g,2-T)}{T},
\end{equation}
with $T=1+|g|^2$.
  \end{teo}

\begin{proof}
Let $h(\mu)=\ell(\mu)$ where $\ell$ be a real function and $\mu=\langle 1,f\rangle$, in this case we have,
\begin{equation*}
h_{,1}=\ell'(\mu) \mu_{,1} \hspace{0.4cm} and \hspace{0.4cm} h_{,2}=\ell'(\mu) \mu_{,2 }
\end{equation*}
\begin{equation*}
h_{,11}=\ell''(\mu) (\mu_{,1}^2+\mu_{,11}) \hspace{0.4cm} and \hspace{0.4cm} h_{,22}= \ell''(\mu) (\mu_{,2}^2+\mu_{,22}).
\end{equation*}
Substituting the above derivatives in (\ref{lapla1}) we have the following equivalence
\begin{equation*}
\ell(\mu) (\ell''(\mu) (\mu_{,1}^2+\mu_{,2}^2)+\ell'(\mu) (\mu_{,11}+ \mu_{,22}) )-C(\mu) (\ell'^2(\mu) (\mu_{,1}^2+\mu_{,2}^2))=0
\end{equation*}

\begin{equation*}
\Leftrightarrow \ell(\mu)\ell''(\mu) |\nabla_L\mu|^2=C(\mu)\ell'^2(\mu) |\nabla_L\mu|^2.
\end{equation*}
Therefore
\begin{equation*}
C(\mu)=\frac{\ell(\mu)\ell''(\mu)}{\ell'^2(\mu)}.
\end{equation*}
On the other hand, we can write (\ref{1}) as
\begin{equation*}
\begin{split}
   X(u) &=\sum_{j=1}^{2} \frac{h_{,j}}{L_{jj}}N_{,j}+hN=\frac{1}{L_{11} }[h_{,1}N_{,1}+h_{,2}N_{,2}] + hN \\
      &=\frac{T^2}{4|g'|^2}\left[ \ell'(\langle 1,f\rangle)\langle 1,f'\rangle \frac{2}{T^2 }(Tg_{,1}-2g\langle g',g\rangle , -2\langle g',g \rangle )\right] \\
      &+\frac{T^2}{4|g'|^2} \left[\ell'(\langle 1,f\rangle)\langle 1,if'\rangle \frac{2}{T^2 }(Tg_{,2}-2g\langle g,ig'\rangle , -2\langle g,ig' \rangle ) \right] +\ell(\langle 1,f\rangle)\frac{(2g,2-T)}{T}\\
&=\frac{\ell'(\langle 1,f\rangle)}{2|g'|^2} (Tg'\bar{f'}-2g\langle g',gf'\rangle,-2 \langle g',gf'\rangle)+\ell(\langle 1,f\rangle)\frac{(2g,2-T)}{T}.
\end{split}
\end{equation*}
\end{proof}

\begin{remark}
The coefficients of the first and second fundamental forms of $X$  and using (\ref{vs})  are given by:
\begin{equation*}
\begin{split}
E&=\langle X_{,1},X_{,1} \rangle=\langle V_{11}N_{,1}+V_{12}N_{,2} ,V_{11}N_{,1}+ V_{12}N_{,2} \rangle\\
&=V_{11}^2L_{11}+V_{12}^2L_{22}=(V_{11}^2+V_{12}^2)L_{11}\\
&= \frac{T^2}{4|g'|^2}[\ell''^2(\mu) \left(\langle 1,f'\rangle^4+\left\langle 1,\frac{if'^2}{2}\right\rangle^2\right)-2\ell''\ell'\left(\langle 1,f'\rangle^2\langle 1,\xi \rangle- \left\langle 1,\frac{if'}{2}\right\rangle\langle i,\xi \rangle\right)\\
& + \ell'^2 (\mu)|\xi|^2]+2\ell[\ell''(\mu) \langle 1,f'\rangle^2- \ell'(\mu)\langle 1,\xi \rangle]+\frac{4\ell^2(\mu)|g'|^2}{T^2},
\end{split}
\end{equation*}

\begin{equation*}
\begin{split}
F&=\langle X_{,1},X_{,2} \rangle=\langle V_{11}N_{,1}+V_{12}N_{,2} ,V_{21}N_{,1}+ V_{22}N_{,2} \rangle\\
&=V_{11}V_{21}L_{11}+V_{12}V_{22}L_{22}=(V_{11}+V_{22}) V_{12}L_{11}\\
&=\left[\frac{\ell''(\mu)|f'|^2T^2}{4|g'|^2} + 2\ell(\mu)] \right] \left[\ell''(\mu) \left\langle 1,\frac{if'^2}{2}\right\rangle + \ell'(\mu)\langle i,\xi \rangle\right],
\end{split}
\end{equation*}

\begin{equation*}
\begin{split}
G&=\langle X_{,2},X_{,2} \rangle=\langle V_{21}N_{,1}+V_{22}N_{,2} ,V_{21}N_{,1}+ V_{22}N_{,2} \rangle=V_{21}^2L_{11}+V_{22}^2L_{22}\\
&=V_{22}^2L_{22}+V_{21}^2L_{11}=(V_{22}^2+V_{21}^2)L_{11}\\
&=\frac{T^2}{4|g'|^2}[\ell''^2(\mu) \left(\langle 1,if'\rangle^4+\left\langle 1,\frac{if'^2}{2}\right\rangle^2\right)+2\ell''\ell'\left(\langle 1,if'\rangle^2\langle 1,\xi \rangle+ \left\langle 1,\frac{if'^2}{2}\right\rangle\langle i,\xi \rangle\right)\\
& +\ell'^2(\mu)|\xi|^2]+2\ell[\ell''(\mu) \langle 1,if'\rangle^2+ \ell'(\mu)\langle 1,\xi \rangle]+\frac{4\ell^2(\mu)|g'|^2}{T^2},
\end{split}
\end{equation*}

\begin{equation*}
\begin{split}
e&=\langle X_{,1},N_{,1} \rangle=\langle V_{11}N_{,1}+V_{12}N_{,2} ,N_{,1}\rangle=V_{ 11}L_{11}\\
&=\ell''(\mu) \langle 1,f'\rangle^2- \ell'(\mu)\langle 1,\xi \rangle + \frac{4\ell(\mu)|g' |^2}{T^2},
\end{split}
\end{equation*}
\begin{equation*}
\begin{split}
f&=\langle X_{,1},N_{,2} \rangle=\langle V_{11}N_{,1}+V_{12}N_{,2} ,N_{,2}\rangle=V_{ 12}L_{22}\\
&=\ell''(\mu) \left\langle 1,\frac{if'^2}{2}\right\rangle + \ell'(\mu)\langle i,\xi \rangle,
\end{split}
\end{equation*}
\begin{equation*}
\begin{split}
g&=\langle X_{,2},N_{,2} \rangle=\langle V_{21}N_{,1}+V_{22}N_{,2} ,N_{,2}\rangle=V_{ 22}L_{22}\\
&= \ell''(\mu) \langle 1,if'\rangle^2+ \ell'(\mu)\langle 1,\xi \rangle+\frac{4\ell(\mu)|g'| ^2}{T^2}.
\end{split}
\end{equation*}
\end{remark}

The following examples are GRT-surfaces for some holomorphic functions $f$ and $g$ with $z=u_1+iu_2\in \mathbb{C}$, using (\ref{trg}) we have

\begin{figure}[htbp]
\begin{center}
\subfigure{\includegraphics[width=29mm]{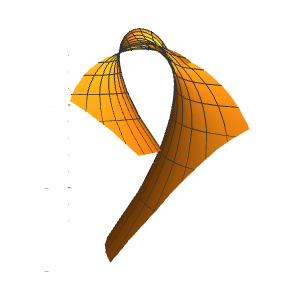}}
\subfigure{\includegraphics[width=29mm]{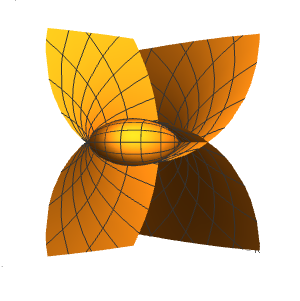}}
\subfigure{\includegraphics[width=29mm]{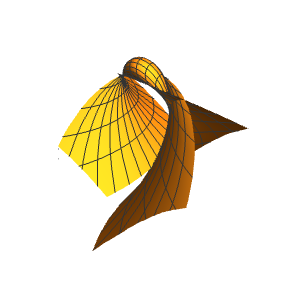}}
\end{center}
\caption{$\ell(t)=t^2+t+1$, \quad $f(z)=g(z)=z=u_1+iu_2$}
\label{lash3}
\end{figure}

\begin{figure}[htbp]
\begin{center}
\subfigure{\includegraphics[width=24mm]{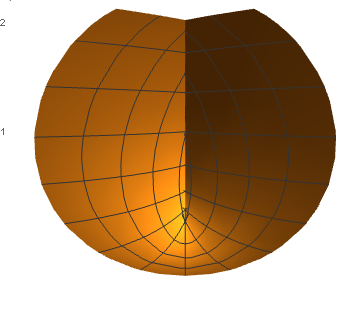}}
\subfigure{\includegraphics[width=24mm]{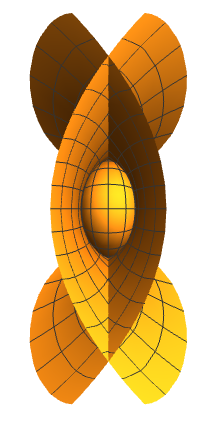}}
\subfigure{\includegraphics[width=24mm]{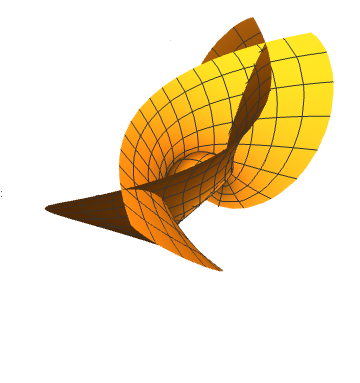}}
\end{center}
\caption{$\ell(t)=\cos(t)$, \quad $f(z)=g(z)=z=u_1+iu_2$}
\label{lash3}
\end{figure}

The following results will be used to classify the GRT-surfaces given by (\ref{trg}).

\begin{lema}
\label{rotation}
Let $\Sigma \subset \mathbb{R}^{n+1} $ be an orientable hypersurface with non-zero Gauss-Kronecker curvature. Then there exists a differentiable function $h:U \rightarrow \mathbb{R}$ such that the normal Gaussian map is given by (\ref{N1}) and $\Sigma $ is locally parameterized by (\ref{1}). Then $X(U)$ is rotational if and only if h is a radial function.
\end{lema}

\begin{proof}
If $X(U)$ is rotational, without loss of generality, we can assume that the axis of rotation is the $x_{n+1}$ axis. This way, the sections orthogonal to the axis of rotation determine in $X(U)$ spheres $(n-1)$ dimensional centered on this axis. Note that along these spheres both $|X|^2$ and the angle between $X$ and $N$ are constant. Given that
\begin{equation*}
\langle X, N \rangle=h, \quad \langle X, X \rangle= |\nabla_Lh|^2 + h^2, \quad L_{ij}=\langle N_{,i}, N_{,j } \rangle.
\end{equation*}
We conclude that $h$ and $|\nabla_Lh|^2$ are constant along these spheres. Taking $N$ as the inverse of the stereographic projection we get,
\begin{equation*}
L_{ij}=\frac{4 \delta _{ij} }{(1+|u|^ 2)^ 2},
\end{equation*}
  so that,
  \begin{equation*}
|\nabla_Lh|^2=\left(\frac{1+|u|^ 2}{2}\right)^2 |\nabla h|,
\end{equation*}
which says that $|u|^2$ is constant as one traverses the orthogonal sections, therefore $h$ is constant along the spheres $(n-1)$ dimensional centered at the origin, therefore $h$ is a function radial.\\

Let $h$ be a radial function, we write $h(u)=J(|u|^2)$, $u\in U $ for $J$ a differentiable function. Let $|u|^2=t$ and denote the derivative of $J$ with respect to $t$ as $J'(t)$. Thus, $h_{,i}=2J'u_i$ and taking $N$ as the inverse of the stereographic projection we have,

\begin{equation*}
X=\left((J'(1-t)+\frac{2J}{1+t})u, -2tJ'+J(\frac{1-t}{1+t}) \right).
\end{equation*}
  If $-2tJ'+J(\frac{1-t}{1+t})$ is constant then
\begin{equation*}
\left|(J'(1-t)+\frac{2J}{1+t})u\right|^2=\left((J'(1-t)+\frac{2J}{1+ t})\right)^2 t.
\end{equation*}
  Which means that the sections orthogonal to the axis $x_{n+1}$ determine in $X(U)$ spheres $(n-1)$ dimensional centered on this axis, so that $X(U)$ is of rotation .
\end{proof}

\begin{teo}
Let $\Sigma$ be a GRT-surface given by (\ref{trg}), connected, where $\ell$ is a real function and $\mu=\langle 1,f\rangle$. Then $\Sigma$ is rotational if and only if, there are constants $a,b\in \mathbb{R}$, such that $\Sigma$ can be locally parameterized by
\begin{equation*}
X_{ab}(u_1,u_2)=\left(M(u_1)\cos(u_2),M(u_1)\sin(u_2),N(u_1)\right),
\end{equation*}
\begin{equation*}
M(u_1)=\frac{1}{1+e^{2u_1}} \left[\frac{a\ell'(\mu)(e^{-u_1}-e^{3u_1})+4\ ell(\mu)e^{u_1}}{2}\right]
\end{equation*}
and
\begin{equation*}
N(u_1)=\frac{1}{1+e^{2u_1}} \left[\ell(\mu)(1-e^{2u_1})-a\ell'(\mu)(1+e ^{2u_1})\right].
\end{equation*}
\end{teo}

\begin{proof}
From Theorem \ref{ttrg}, we have that $\Sigma$ is locally parameterized by
\begin{equation*}
X=\frac{\ell'(\mu)}{2|g'|^2}(Tg'\bar{f'}-2g\langle g',gf'\rangle , -2\langle g', gf'\rangle) + \ell(\mu)\frac{(2g,2-T)}{T},
\end{equation*}
where $f,g$ are holomorphic functions, $\ell$ is not constant and remembering that $h=\ell(\mu)$, $h,_2=\ell'(\mu)\langle 1,if'\rangle=0$ and $\langle 1,if'\rangle=0$. The Cauchy-Riemann equations guarantee us that $f(z)=az+z_0$, with $z=u_1+iu_2, z_0=b+ic\in \mathbb{C}$, thus Thus, by the lemma \ref{rotation} we have that $\Sigma$ is rotational if and only if
\begin{equation*}
g(z)=e^z, \hspace{.2cm} h(z)=\ell(au_1+b),
\end{equation*}
so that $f'(z)=a$, $g'(z)=g(z)$ and $T=1+e^{2u_1}$. In this conditions,
\begin{equation*}
\begin{split}
    X &=\frac{\ell'(\mu)}{2|g'|^2}(Tg'\bar{f'}-2g\langle g',gf'\rangle , -2\langle g' ,gf'\rangle) + \ell(\mu)\frac{(2g,2-T)}{T} \\
      &=\frac{\ell'(\mu)}{2e^{u_1}}\left((1+e^{2u_1})e^za-2e^z\langle e^z,e^za\rangle , -2\langle e^z , e^za\rangle \right)+\ell(\mu)\left(\frac{2e^z}{1+e^{2u_1}},\frac{1-e ^{2u_1}}{1+e^{2u_1}}\right)\\
      &=\frac{\ell'(\mu)}{2e^{2u_1}}\left((a+ae^{2u_1})e^z-2ae^{2u_1}e^z,-2ae^{2u_1 } \right)+\ell(\mu)\left(\frac{2e^z}{1+e^{2u_1}},\frac{1-e^{2u_1}}{1+e^{2u_1} }\right)\\
      &=\frac{\ell'(\mu)}{2e^{2u_1}}\left((a-ae^{2u_1})e^z,-2ae^{2u_1} \right)+\ell(\mu)\left(\frac{2e^z}{1+e^{2u_1}},\frac{1-e^{2u_1}}{1+e^{2u_1}}\right)\\
      &=\ell'(\mu)\left(\frac{ae^z(e^{-2u_1}-1)}{2},-a\right)+\ell(\mu)\left(\frac {2e^z}{1+e^{2u_1}},\frac{1-e^{2u_1}}{1+e^{2u_1}}\right)\\
      &=\left(\frac{\ell'(\mu)ae^z(e^{-2u_1}-1)}{2},-\ell'(\mu)a\right)+\left(\frac{2\ell(\mu)e^z}{1+e^{2u_1}},\frac{\ell(\mu)(1-e^{2u_1})}{1+e^{2u_1} }\right)\\
      &=\left(\left(\frac{a\ell'(\mu)(e^{-u_1}-e^{u_1})}{2}+\frac{2\ell(\mu)e^ {u_1}}{1+e^{2u_1}}\right)(\cos(u_2)+i\sin(u_2)),\frac{\ell(\mu)(1-e^{2u_1})} {1+e^{2u_1}}-a\ell'(\mu)\right).
\end{split}
\end{equation*}
follow or result.
\end{proof}

The following figures are examples of rotating GRT-surfaces for some $\ell$ functions and constants $a$ and $b$.

\begin{figure}[htbp]
\begin{center}
\subfigure{\includegraphics[width=40mm]{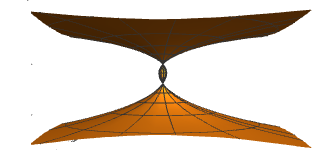}}
\subfigure{\includegraphics[width=40mm]{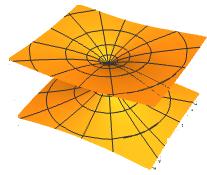}}
\end{center}
\caption{$\ell(t)=t^2+t+1$, \quad $a=1,\hspace{0.2cm} b=0$}
\label{lash3}
\end{figure}

\begin{figure}[htbp]
\begin{center}
\subfigure{\includegraphics[width=40mm]{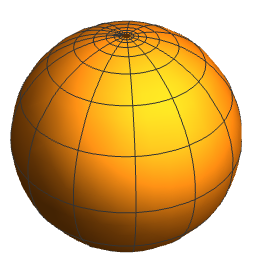}}
\end{center}
\caption{$\ell(t)=t^2+t+1$, \quad $a=0,\hspace{0.2cm} b=1$}
\label{lash3}
\end{figure}

\begin{figure}[htbp]
\begin{center}
\subfigure{\includegraphics[width=40mm]{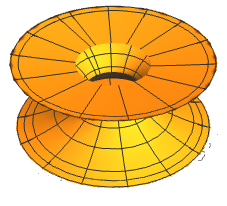}}
\subfigure{\includegraphics[width=40mm]{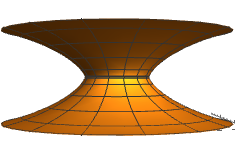}}
\end{center}
\caption{$\ell(t)=\cos(t)$, \quad $a=1,\hspace{0.2cm} b=0$}
\label{lash3}
\end{figure}


\begin{figure}[htbp]
\begin{center}
\subfigure{\includegraphics[width=40mm]{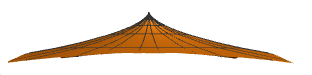}}
\subfigure{\includegraphics[width=40mm]{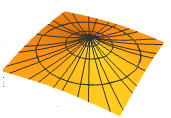}}
\end{center}
\caption{$\ell(t)=\sinh(t)$, \quad $a=1,\hspace{0.2cm} b=0$}
\label{lash3}
\end{figure}


\newpage
\renewcommand{\refname}{Bibliografía}








\end{document}